\documentclass[a4paper,11pt]{amsart}
\usepackage[english]{babel}
\usepackage{amstext,amsfonts,amsthm,indentfirst,graphicx,amssymb,amscd,epsfig}
\usepackage{amsmath}
\usepackage{graphics,wrapfig}
\usepackage[ansinew]{inputenc}    
\usepackage{bm}
\usepackage{mathrsfs,dsfont}
\usepackage{hyperref}

\newtheorem{ex}{Example}[section]
\newtheorem{proposition}{Proposition}[section]
\newtheorem{definition}{Definition}[section]
\newtheorem{theorem}{Theorem}[section]
\newtheorem{corollary}{Corollary}[section]
\newtheorem{remark}{Remark}[section]
\newtheorem{lemma}{Lemma}[section]

\newfont{\bbb}{msbm10 scaled\magstephalf}     

\def\R{\mathbb R}

\def\R{\mbox{\bbb R}}
\def\x{\mathbf{x}}
\def\y{\mathbf{y}}

\def\EO{E_\Omega}
\def\FO{F_\Omega}
\def\GO{G_\Omega}
\def\eo{e_\Omega}
\def\fuo{f_{1\Omega}}
\def\fdo{f_{2\Omega}}
\def\go{g_\Omega}
\def\n{\mathbf{n}}
\def\p{\mathbf{p}}
\def\q{\mathbf{q}}
\def\I{\mathbf{I}}

\def\II{\mathbf{II}}

\def\w{\mathbf{w}}
\def\id{\mathbb{I}}

\def\Omegam{\mathbf{\Omega}}

\def\Lambdam{\mathbf{\Lambda}}

\def\mum{\boldsymbol{\mu}}

\newcommand{\spann}[2]{\left\langle{#1},{#2}\right\rangle}

\newcommand{\laO}{\lambda_{\Omega}}
\newcommand{\KO}{K_{\Omega}}
\newcommand{\HO}{H_{\Omega}}

\title[Extendibility and boundedness of invariants on wavefronts]{Extendibility and boundedness of invariants on singularities of wavefronts}

\author{T. A. Medina-Tejeda}

\date{}

\address{Instituto de Ci\^encias Matem\'aticas e de Computa\c{c}\~ao - Universidade de S\~ao Paulo,
Av. Trabalhador s\~ao-carlense, 400 - Centro,
CEP: 13566-590 - S\~ao Carlos - SP, Brazil}

\email{tamedinat@usp.br}

\thanks{The author is supported by CAPES Grant no. PROEX-10359340/D}

\begin{document}
\begin{abstract}
We investigate necessary and sufficient conditions for the extendibility and boundedness of Gaussian curvature, Mean curvature and principal curvatures near all types of singularities on fronts. We also study the convergence to infinite limits of these geometrical invariants and show how this is tightly related to a particular property of uniform approximation of fronts by parallel surfaces.    
\end{abstract}

\maketitle

\section{Introduction}
Frontal is a class of surfaces with singularities that have been intensely studied in these last years within the field of differential geometry and singularities, in particular wave fronts \cite{agv,arnld-sing-caus,ishifrontal2,MSUY,med,murataumehara,front,t1,t2,t3}. They admit a normal vector field along their parametrizations defined even on singularities of these. More precisely, a smooth map $\x: U \to \R^3$ defined in an open set $U \subset \R^2$ is called a {\it frontal} if, for all $\p\in U$ there exists a unit normal vector field $\n_p:V_p \to \R^3$ along $\x$ (i.e $\x_u$, $\x_v$ are orthogonal to $\n$), where $V_p$ is an open set of $U$, $\p\in V_p$. If the {\it singular set} $\Sigma(\x)=\{\p\in U: \x \text{ is not immersive at $\p$} \}$ has empty interior we call $\x$ a {\it proper frontal} and if $(\x,\n_p):U \to \R^3\times S^2$ is an immersion for all $\p\in U$ we call $\x$ a {\it wave front} or simply {\it front}.       
 
The behavior of Gaussian curvature, Mean curvature and principal curvatures near non-degenerate singularities on wave fronts have been widely studied in \cite{front,MSUY,t1,t3}. However, in the degenerate case this is unknown, as well as the convergence to infinite limits of these invariants has been little explored. For this reason it is natural to wonder which properties of wave fronts determine one behavior or another on general types of singularities. Also, there is a lack of literature about the geometry of singularities of rank 0 (or corank 2) on wave fronts and our approach here allow us to study them. For a singularity $\p$ of a wavefront $\x$, there exists $l>0$ and a neighborhood $U_l$ of $\p$ such $\y_l=\x+l\n$ is an immersion, also this $l$ can be chosen as small as we wish. The neighborhood $U_l$ may shrink as $l$ is smaller, then is natural to ask when $U_l$ can be hold fixed for $l$ arbitrarily small, in this case we say that $\x$ is {\it parallelly smoothable at $\p$}. We will see that, this last property is determined by the convergence to infinite limits of the geometric invariants at each type of singularity and also is related with the extendibility of the principal curvatures at singularities of rank 1.  
 
In section 2, we establish the notation, terminology and basic results that we use. Most of it is introduced and obtained in \cite{med}. In section 3 we introduce the relatives principal curvatures which give us geometrical information near singularities and are defined even on them. In section 4, we study singularities of rank 1 both degenerate and non-degenerate. The theorems \ref{K1}, \ref{K2} give equivalent conditions for boundedness and extendibility of the Gaussian curvature at degenerate and non-degenerate singularities which generalize the one found in \cite{front} and similarly theorem \ref{ka}  for the principal curvatures. By last in section 5, we study the behavior of the invariants at singularities of rank 0, in particular the boundedness and extendibility of the Mean curvature obtaining results quite different from those obtained in the rank 1 case. The theorems \ref{pa1}, \ref{pa2} and \ref{p2} characterize when a wavefront is parallelly smoothable at all types of singularities, which is determined by the convergence to infinite limits of the geometric invariants. The example \ref{e1} shows explicitly a wavefront with singularity of rank $0$ and extendable Mean curvature.

\section{Fixing notation, definitions and some basic results}\label{section-definition}
In this paper, all the maps and functions are of class $C^\infty$. We denote $U$ and $V$ open sets in $\R^2$ when is not mentioned anything about them. Let $\x:U \to \R^3$ be a smooth map, we call a {\it tangent moving basis} (tmb) of $\x$ a smooth map $\Omegam:U \to \mathcal{M}_{3\times2}(\R)$ in which the columns $\w_1, \w_2:U \to \R^3$ of the matrix $\Omegam=\begin{pmatrix}\w_1&\w_2\end{pmatrix}$ are linearly independent smooth vector fields and  $\x_u,\x_v \in \spann{\w_1}{\w_2}$, where $\spann{}{}$ denotes the linear span vector space. It is known that a smooth map $\x:U \to \R^3$ is a frontal if and only if there exist tangent moving bases of $\x$ locally. Since we are interested in exploring local properties of frontals, we always assume that we have a global tmb $\Omegam$ for $\x$. We denote by $\n:=\frac{\w_1\times\w_2}{\|\w_1\times\w_2\|}$ the normal vector field induced by $\Omegam$. Let $\mathbf{f}:U \to \R^n$ be a smooth map, we denote by $D\mathbf{f}:=(\frac{\partial \mathbf{f}_i}{\partial x_j})$, the differential of $\mathbf{f}$ and we consider it as a smooth map $D\mathbf{f}:U \to \mathcal{M}_{n\times2}(\R)$. We write $D\mathbf{f}_{x_1}$, $D\mathbf{f}_{x_2}$ the partial derivatives of $D\mathbf{f}$ and $D\mathbf{f}(\p):=(\frac{\partial \mathbf{f}_i}{\partial x_j}(\p))$ for $\p \in U$. Also, vectors in $\R^n$ are identified as column vectors in $\mathcal{M}_{n\times1}(\R)$ and if $\mathbf{A}\in \mathcal{M}_{n\times n}(\R)$, $\mathbf{A}_{(i)}$ is the $i^{th}$-row and $\mathbf{A}^{(j)}$ is the $j^{th}$-column of $\mathbf{A}$. The trace and adjoint of a matrix are denoted by $tr()$ and $adj()$ respectively. 

Let $\x:U \to \R^3$ be a frontal, $\Omegam$ a tmb of $\x$. Denoting $()^T$ the operation of transposing a matrix, we set the matrices of the first and second fundamental forms: 
$$\mathbf{I}=
\begin{pmatrix}
E&F\\
F&G
\end{pmatrix}:=D\x^TD\x,\ \mathbf{II}=\begin{pmatrix}
e&f\\
f&g
\end{pmatrix}:=-D\x^T D\n.$$	
The Weingarten matrix $\boldsymbol{\alpha}:=-\II^{T}\I^{-1}$ is defined in $\Sigma(\x)^c$. Also, we set the matrices:
$$\I_\Omega=\begin{pmatrix}
\EO & \FO\\
\FO & \GO\end{pmatrix}:=\Omegam^T\Omegam\label{IO},\ \II_\Omega=\begin{pmatrix}
\eo&\fuo\\
\fdo&\go
\end{pmatrix}:=-\Omegam^T D\n\label{IIO},$$
$$\boldsymbol{\mu}_{\Omega}:=-\II_\Omega^T\I_\Omega^{-1}\label{W},\ \Lambdam_{\Omega}:=D\x^T\Omegam(\I_\Omega)^{-1},\ \boldsymbol{\alpha}_\Omega:=\boldsymbol{\mu}_{\Omega}adj(\Lambdam_\Omega).$$
If $\x$ is a frontal and $\Omegam$ a tmb of $\x$, we write simply $\Lambdam=(\lambda_{ij})$ and $\mum=(\mu_{ij})$ instead of $\Lambdam_\Omega$ and $\boldsymbol{\mu}_{\Omega}$ when there is no risk of confusion, $\laO:=det(\Lambdam)$ and $\mathfrak{T}_\Omega(U)$ as the principal ideal generated by $\laO$ in the ring $C^\infty(U,\R)$. The matrix $\Lambdam$ and $\boldsymbol{\mu}$ satisfy $D\x=\Omegam\Lambdam^T$ and $D\n=\Omegam\mum^T$ (see \cite{med}), thus $\Sigma(\x)=\laO^{-1}(0)$ and $rank(D\x)=rank(\Lambdam)$.

\begin{definition}
	Let $\x:U \to \R^3$ be a frontal, $\Omegam$ a tmb of $\x$, $\p \in \Sigma(\x)$, we say that $\p$ is a {\it non-degenerate} singularity if $D\laO(\p)\neq(0,0)$, in another case is called degenerate.
\end{definition}
\begin{remark}
	This definition does not depend on the chosen tmb $\Omegam$ because with another tmb $\bar{\Omegam}$, we have $\lambda_{\bar{\Omega}}=\rho\laO$ with $\rho:U \to \R$ being a $C^{\infty}$-function that never vanish.
\end{remark} 

The following propositions were proved in \cite{med} and we are going to use these frequently.
\begin{definition}
	Let $\x:U \to \R^3$ be a frontal, $\Omegam$ a tmb of $\x$, the {\it $\Omega$-relative curvature} and the {\it $\Omega$-relative mean curvature} are defined on $U$ by $\KO=det(\boldsymbol{\mu}_{\Omega})$ and $H_\Omega=-\frac{1}{2}tr(\boldsymbol{\alpha}_\Omega)$ respectively.
\end{definition}

\begin{proposition}\cite[Proposition 3.17]{med}\label{lim} 
	Let $\x:U \to \R^3$ be a proper frontal, $\Omegam$ a tangent moving basis  of $\x$, $\KO$, $H_\Omega$, $K$ and $H$ the $\Omega$-relative curvature, the $\Omega$-relative mean curvature, the Gaussian curvature and the mean curvature of $\x$ respectively.Then,	
	\begin{enumerate}
		\item for $\p\in \Sigma(\x)^c$, $\KO=\laO K$ and $H_\Omega=\laO H$,
		\item for $\p\in \Sigma(\x)$, $\KO=\lim\limits_{(u,v)\to p}\laO K$ and $H_\Omega=\lim\limits_{(u,v)\to p}\laO H$, 
		
	\end{enumerate} 
	where the right sides are restricted to the open set $\Sigma(\x)^c$.	
\end{proposition}

\begin{proposition}\cite[Theorem 3.5]{med}\label{D}
	Let $\x:U \to \R^3$ be a frontal and $\Omegam$ a tmb of $\x$, then first and second fundamental forms have the following decomposition:
$$\I=\Lambdam\I_{\Omega}\Lambdam^T,\ \II=\Lambdam\II_{\Omega}$$		
which satisfy,
	\begin{subequations}
		\begin{align*}
		\Lambdam_{(1)u}^{\phantomsection}\I_\Omega\Lambdam_{(2)}^T-\Lambdam_{(1)}^{\phantomsection}\I_\Omega\Lambdam_{(2)u}^T+E_v-F_u \in \mathfrak{T}_\Omega\\
		\Lambdam_{(1)v}^{\phantomsection}\I_\Omega\Lambdam_{(2)}^T-\Lambdam_{(1)}^{\phantomsection}\I_\Omega\Lambdam_{(2)v}^T+F_v-G_u \in \mathfrak{T}_\Omega
		\end{align*}
	\end{subequations}
\end{proposition}

\begin{proposition}\cite[Theorem 3.22]{med}\label{wft}
	Let $\x:U \to \R^3$ be a frontal, $\Omegam$ a tangent moving basis of $\x$ and  $\p \in \Sigma(\x)$. Then,
	\begin{enumerate}
		\item $\x:U \to \R^3$ is a front on a neighborhood $V$ of $\p$ with $rank(D\x(\p))=1$ if and only if $H_\Omega(\p)\neq 0$. 
		
		\item $\x:U \to \R^3$ is a front on a neighborhood $V$ of $\p$ with $rank(D\x(\p))=0$ if and only if $H_\Omega(\p)=0$ and $\KO(\p)\neq0$.
		
	\end{enumerate}
\end{proposition}

\begin{corollary}\label{wf}
	A frontal $\x:U \to \R^3$ is a front if and only if $(\KO,H_\Omega)\neq(0,0)$ on $\Sigma(\x)$ for whatever tangent moving basis $\Omegam$ of $\x$.
\end{corollary} 

\section{The relative principal curvatures}

\begin{proposition}\label{Weinr} 
	Let $\x:U \to \R^3$ be a proper frontal, $\Omegam$ a tmb of $\x$, we have the following equality, $\boldsymbol{\alpha}_\Omega=\boldsymbol{\alpha}\laO$ on $\Sigma(\x)^c$, in particular $\boldsymbol{\alpha}_\Omega$ has real eigenvalues.
\end{proposition}
\begin{proof}
	By proposition \ref{D}, $\I=\Lambdam\I_\Omega \Lambdam^T$ and $\II=\Lambdam\II_\Omega$, then for $\p\in \Sigma(\x)^c$, $\boldsymbol{\alpha}=-\II^T\I^{-1}=-\II_\Omega^T\Lambdam^T(\Lambdam^T)^{-1}\I_\Omega^{-1}\Lambdam^{-1}=\boldsymbol{\mu}_{\Omega}\Lambdam^{-1}$. Thus, we have $\boldsymbol{\alpha}\laO=\boldsymbol{\mu}_{\Omega}adj(\Lambdam)=\boldsymbol{\alpha}_\Omega$. The eigenvalues of $\boldsymbol{\alpha}_\Omega$ are real if $tr(\boldsymbol{\alpha}_\Omega)^2-4 det(\boldsymbol{\alpha}_\Omega)\geq0$. As $\KO=det(\boldsymbol{\mu}_{\Omega})$ and $H_\Omega=-\frac{1}{2}tr(\boldsymbol{\alpha}_\Omega)$ this is equivalent to have $\HO^2-\laO\KO\geq0$ and since $\HO^2-\laO\KO=\laO^2(H^2-K)\geq0$ on $\Sigma(\x)^c$, then by continuity and the density of regular points, it follows the result on $U$.      
\end{proof}
Denoting the eigenvalues of $\boldsymbol{\alpha}_{\Omega}$ by $-k_{1\Omega}$, $-k_{2\Omega}$, then $k_{1\Omega}$, $k_{2\Omega}$ satisfy the equation $k^2+tr(\boldsymbol{\alpha}_\Omega^T) k+det(\boldsymbol{\alpha}_\Omega^T)=0$. Since $\KO=det(\boldsymbol{\mu}_{\Omega})$ and $H_\Omega=-\frac{1}{2}tr(\boldsymbol{\alpha}_\Omega)$, we have $k^2-2H_\Omega k+\KO=0$. Thus,
$$k=\HO\pm\sqrt{\HO^2-\laO\KO}$$

\begin{definition}\label{pk}
Let $\x:U \to \R^3$ be a frontal, $\Omegam$ a tmb of $\x$, we call the functions $k_{1\Omega}:=\HO-\sqrt{\HO^2-\laO\KO}$ and $k_{2\Omega}:=\HO+\sqrt{\HO^2-\laO\KO}$ the {\it relative principal curvatures}. We also define the following functions on $U-\Sigma(\x)$:
\begin{equation*}
k_1:= \begin{cases}
H-\sqrt{H^2-K} & \text{if } \laO>0,\\
H+\sqrt{H^2-K} & \text{if } \laO<0.
\end{cases}
\end{equation*}	
\begin{equation*}
k_2:= \begin{cases}
H+\sqrt{H^2-K} & \text{if } \laO>0,\\
H-\sqrt{H^2-K} & \text{if } \laO<0.
\end{cases}
\end{equation*}
\end{definition}  
We clarify that, the principal curvatures of $\x$ are the functions defined by $\kappa_{-}:=H-\sqrt{H^2-K}$ and $\kappa_{+}:=H+\sqrt{H^2-K}$ on $U-\Sigma(\x)$.

\begin{remark}
	$k_1$ and $k_2$ are smooth functions on $U-\Sigma(\x)$ and they have similar properties to the classical principal curvatures. Also their definitions do not depend on the chosen tmb $\Omegam$ inducing the same orientation of the normal vector field $\n$. If another tmb $\hat{\Omegam}$ induces an opposite orientation of $\n$, then the signs of these functions are opposite as well. Observe that $k_1k_2=K$ and $\frac{k_1+k_2}{2}=H$ on $U-\Sigma(\x)$. In the case of non-degenerate singularities, if we make a suitable change of coordinates $k_1, k_2$ coincides with those functions defined in (\cite{t1}, equation (2.6)). 
\end{remark}

\begin{proposition}\label{lim2} 
	Let $\x:U \to \R^3$ be a proper frontal, $\Omegam$ a tangent moving basis of $\x$. Then,	
	\begin{enumerate}
		\item for $\p\in \Sigma(\x)^c$, $k_{1\Omega}=\laO k_1$ and $k_{2\Omega}=\laO k_2$,
		\item for $\p\in \Sigma(\x)$, $k_{1\Omega}=\lim\limits_{(u,v)\to p}\laO k_1$ and $k_{2\Omega}=\lim\limits_{(u,v)\to p}\laO k_2$. 
		
	\end{enumerate} 
		
\end{proposition}

\begin{proof}
	We have that
	$k_{1\Omega}=\laO H-\sqrt{\laO^2 H^2-\laO^2 K}=\laO H-|\laO|\sqrt{ H^2- K}=\laO k_1$ and similarly $k_{2\Omega}=\laO k_2$ on $\Sigma(\x)^c$. For $\p\in \Sigma(\x)$, by smoothness of $k_{1\Omega}$, $k_{2\Omega}$ and density of $\Sigma(\x)^c$, $k_{1\Omega}=\lim\limits_{(u,v)\to p}\laO k_1$ and $k_{2\Omega}=\lim\limits_{(u,v)\to p}\laO k_2$. 
\end{proof}

In \cite{med}(proof of proposition 3.19) was observed that making change of coordinates $\mathbf{h}$ on a frontal $\x$ and taking $\hat{\Omegam}:=\Omegam \circ \mathbf{h}$ as tmb of $\x \circ \mathbf{h}$, it results with new different relative curvatures $det(\mathbf{h})(\KO \circ \mathbf{h})$ and $det(\mathbf{h})(\HO \circ \mathbf{h})$. However, if we choose the tmb $\Omegam^h:=(\Omegam \circ \mathbf{h})D\mathbf{h}$ instead of $\hat{\Omegam}$ when we make a change of coordenates, they remain invariant. 

\begin{proposition}[Invariance property]\label{inv}
	Let $\x:U \to \R^3$ be a frontal, $\Omegam$ a tmb of $\x$, $\mathbf{h}:V \to U$ diffeomorphism, then the new relative curvatures of $\x \circ \mathbf{h}$ are $K_{\Omega^h}=\KO\circ\mathbf{h}$ and $H_{\Omega^h}=\HO\circ\mathbf{h}$. In particular, $k_{1\Omega^h}=k_{1\Omega}\circ \mathbf{h}$, $k_{2\Omega^h}=k_{2\Omega}\circ \mathbf{h}$ and $\lambda_{\Omega^h}=\laO\circ\mathbf{h}$.   
\end{proposition}
\begin{proof}
	Observe that, the matrix $\Lambdam_{\Omega^h}$ induced by $\Omegam^h$ is $(D\mathbf{h})^{-1}(\Lambdam_{\Omega}\circ \mathbf{h}) D\mathbf{h}$, then $\lambda_{\Omega^h}=\laO\circ\mathbf{h}$ and 
	since $D(\n\circ\mathbf{h})=\Omegam^h\boldsymbol{\mu}_{\Omega^h}^T$, we have $(D\mathbf{h})^{-1}\boldsymbol{\mu}_{\Omega}^T(D\mathbf{h})=\boldsymbol{\mu}_{\Omega^h}^T$ and therefore $K_{\Omega^h}=\KO\circ\mathbf{h}$. Also, $\boldsymbol{\mu}_{\Omega^h}adj(\Lambdam_{\Omega^h})=det(D\mathbf{h})D\mathbf{h}^{-1}\boldsymbol{\mu}_{\Omega}adj(\Lambdam_{\Omega})adj(D\mathbf{h}^{-1})$ and using that $tr(\mathbf{A}\mathbf{B}\mathbf{C})=tr(\mathbf{C}\mathbf{A}\mathbf{B})$ we get $H_{\Omega^h}=\HO\circ\mathbf{h}$. Since, $k_{1\Omega^h}$ and $k_{2\Omega^h}$ are written in terms of $K_{\Omega^h}$ and $H_{\Omega^h}$, the result follows for them.      
\end{proof}     

\begin{lemma}\label{lem1}
	Let $\x:U \to \R^3$ be a wavefront, $\Omegam$ a tmb of $\x$, for each $\p \in \Sigma(\x)$ there exist locally an embedding $\y_l:V \to \R^3$, $\p \in V$, such that $D\y_l$ is a tmb of $\x$, the matrix $\Lambdam_{D\y_l}$ determined for this tmb is $\id-l\boldsymbol{\alpha}_l$, where $l\in \R^+$ and $\boldsymbol{\alpha}_l$ is the Weingarden matrix of $\y_l$ and $\id$ is the identity matrix. 
\end{lemma}
\begin{proof}
	For each $t\in \R$, consider $\y_t=\x+t\n$, as $D\x=\Omegam\Lambdam^T$ and $D\n=\Omegam\mum^T$ we have $D\y_t=\Omegam\Lambdam^T+t\Omegam\boldsymbol{\mu}^T$, then $\y_t$ has a singularity at $\q$ if and only if $det(\Lambdam^T+t\boldsymbol{\mu}^T)(\q)=0$. Making a direct computation $det(\Lambdam^T+t\boldsymbol{\mu}^T)=\laO-2t\HO+t^2\KO$ and now taking $\p\in\Sigma(\x)$, by corollary (\ref{wf}), there exist $l\in \R^+$ such that $det(\Lambdam^T+l\boldsymbol{\mu}^T)(\p)=-2l\HO(\p)+l^2\KO(\p)\neq 0$. Thus, there exist a neighborhood $V$ of $\p$ such that $\y_{l}:V\to \R^3$ is an embedding. Since, $D\y_{l}=\Omegam(\Lambdam+l\mum)$, $D\y_{l}$ is a tmb of $\x$. We can assume $D\y_{l}$ and $\Omegam$ induce the same normal vector $\n$ (i.e $det(\Lambdam^T+l\boldsymbol{\mu}^T)>0$ on $V$), otherwise we can change the order of column in $\Omegam$ from the beginning. Therefore, we have $D\y_{l}(\id-l\boldsymbol{\alpha}^T_l)=D\y_{l}-lD\n=D\x=D\y_{l}\Lambdam_{D\y_l}^T$, thus $\Lambdam_{D\y_l}=\id-l\boldsymbol{\alpha}_l$.    
\end{proof}

\begin{lemma}\label{lem2}
	Let $\x:U \to \R^3$ be a wavefront, $\p \in \Sigma(\x)$ and $\Omegam=D\y_l$ a tmb of $\x$ as above, then  
	\begin{enumerate}
		\item $\I_\Omega=\I_l$, $\II_\Omega=\II_l$, $\boldsymbol{\mu}=\boldsymbol{\alpha}_l$ and $\boldsymbol{\mu}adj(\Lambdam)=\boldsymbol{\alpha}_l-lK_l\id$.
		
		\item $\KO=K_l$, $\HO=H_l+K_ll$ and $\laO=1+2H_ll+K_ll^2$.
		
		\item $k_{1\Omega}=k_{1l}(1+lk_{2l})$ and $k_{2\Omega}=k_{2l}(1+lk_{1l})$.
	\end{enumerate}
	 where $\I_l$, $\II_l$, $K_l$, $H_l$, $k_{1l}$, $k_{2l}$ are first fundamental form, second fundamental form, Gaussian curvature, mean curvature and principal curvatures of $\y_l$ respectively. Additionally, $rank(D\x(\p))=1$ if and only if $\y_l$ is free of umbilical point on a neighborhood of $\p$. Similarly, $rank(D\x(\p))=0$ if and only if $\y_l$ has a umbilical point at $\p$ of positive Gaussian curvature. 
\end{lemma}
\begin{proof}
\	
\begin{enumerate}
	
\item Applying the definition directly we get the first three equalities. By lemma \ref{lem1} $\Lambdam=\id-l\boldsymbol{\alpha}_l$, then $\boldsymbol{\mu}adj(\Lambdam)=\boldsymbol{\alpha}_l (\id-l adj(\boldsymbol{\alpha}_l))=\boldsymbol{\alpha}_l-lK_l\id$.  
		
\item Using item $(1)$, $\KO=det(\mum)=det(\boldsymbol{\alpha}_l)=K_l$, $\HO=-\frac{1}{2}tr(\mum adj(\Lambdam))=-\frac{1}{2}tr(\boldsymbol{\alpha}_l-lK_l\id)=H_l+K_ll$ and $\laO=det(\Lambdam)=det(\id-l\boldsymbol{\alpha}_l)=1+2H_ll+K_ll^2$.
		
\item Using the formulas in definition \ref{pk}, item $(2)$ and knowing that $k_{1l}=H_l-\sqrt{H_l^2-K_l}$, $k_{2l}=H_l+\sqrt{H_l^2-K_l}$ and $K_l=k_{1l}k_{2l}$, a simple computation leads to item $(3)$.
\end{enumerate}
For the last part, by proposition \ref{wft} $rank(D\x(\p))=0$ if and only if $\HO(\p)=0$, $\laO(\p)=0$ and $\KO(\p)\neq 0$. On the other hand these conditions are equivalent to $H_l(\p)=-K_l(\p)l$ and $0=1-2l^2K_l(\p)+l^2K_l(\p)$ which is equivalent to $K_l(\p)=\frac{1}{l^2}$ and $H_l(\p)=-\frac{1}{l}$. Then $\y_l$ has an umbilical point at $\p$ of positive Gaussian curvature. Conversely, we have $0<K_l(\p)=H_l(\p)^2$, then $0=1+2lH_l(\p)+l^2H_l(\p)^2$ which imply $H_l(\p)=-\frac{1}{l}$ and therefore $\HO(\p)=-\frac{1}{l}+\frac{1}{l^2}l=0$, by proposition \ref{wft} $rank(D\x(\p))=0$. Equivalently $rank(D\x(\p))=1$ if and only if $\p$ is not a umbilical point which is equivalent to have $\y_l$ free of umbilical point on a neighborhood of $\p$.          
\end{proof}

\section{Singularities of rank 1}

In the following propositions, we study the behavior at a singular point of rank 1 of the classical invariants and the functions defined in previous section. The non-degenerate case was investigated in \cite{front,MSUY,t1,t3}. The results here include degenerate and non-degenerate cases.

\begin{proposition}\label{prin1}
Let $\x:U \to \R^3$ be a proper wavefront, $\Omegam$ a tmb of $\x$, then for every $\p \in \Sigma(\x)$ with $rank(D\x(\p))=1$, the following is always satisfied:
\begin{enumerate}
	\item $(k_{1\Omega}(\p),k_{2\Omega}(\p))\neq(0,0)$. In particular, if $k_{1\Omega}(\p)\neq 0$ (resp. $k_{2\Omega}(\p)\neq 0$), then $k_{2\Omega}(\p)=0$ (resp. $k_{1\Omega}(\p)=0$). Also, $\HO(\p)<0$ (resp. $\HO(\p)>0$) if and only if $k_{1\Omega}(\p)\neq 0$ (resp. $k_{2\Omega}(\p)\neq 0$). 
	
	\item There is an open neighborhood $V \subset U$ of $\p$ in which one of the functions $k_1$, $k_2$ has a $C^\infty$ extension to $V$. More precisely, $k_1$ (resp. $k_2$) has a $C^\infty$ extension if only if $k_{1\Omega}(\p)=0$ (resp. $k_{2\Omega}(\p)$). 
	
	\item One of the functions $k_1$, $k_2$ in absolute value converge to $\infty$. More precisely, $\lim\limits_{(u,v)\to p}|k_1|=\infty$ (resp. $|k_2|$) if and only if $k_{1\Omega}(\p)\neq 0$ (resp. $k_{2\Omega}(\p)$).
	
	\item $\lim\limits_{(u,v)\to p}|H|=\infty$.
	
	\item If $\KO(\p)\neq 0$ then  $\lim\limits_{(u,v)\to p}|K|=\infty$.  
	
\end{enumerate} 
	
\end{proposition}
\begin{proof}
\
\begin{enumerate}

\item Observe that, for $\p \in \Sigma(\x)$, $k_{1\Omega}(\p)=\HO(\p)-|\HO(\p)|$ and $k_{2\Omega}(\p)=\HO(\p)+|\HO(\p)|$. By proposition \ref{wft}, $\HO(\p)\neq 0$, then just one of $k_{1\Omega}(\p)$, $k_{2\Omega}(\p)$ is different of zero. Thus, the sub index of $k_{i\Omega}(\p)$ corresponding to the non-zero value is determined bijectively by the sign of $\HO(\p)$. 

\item By item $(1)$, without loss of generality, we can assume that $k_{1\Omega}(\p)\neq 0$. Let $V$ a neighborhood of $\p$ such that $k_{1\Omega}\neq 0$ on $V$, then by proposition \ref{lim}, $k_2=\frac{\laO k_1k_2}{\laO k_1}=\frac{\KO}{k_{1\Omega}}$ on $V-\Sigma(\x)$. Thus, $\frac{\KO}{k_{1\Omega}}$ is a $C^{\infty}$ extension of $k_2$ to $V$. 

\item By item $(1)$, without loss of generality, we can assume that $k_{1\Omega}(\p)\neq 0$. Let $V$ be a neighborhood of $\p$ such that $k_{1\Omega}\neq 0$ on $V$, then $k_1=\frac{k_{1\Omega}}{\laO}$ on $V-\Sigma(\x)$. Thus, for every $\p \in \Sigma(\x)\cap V$, $\lim\limits_{(u,v)\to p}|k_1|=\lim\limits_{(u,v)\to p}\frac{|k_{1\Omega}|}{|\laO|}=\infty$.  

\item Since $\HO(\p)\neq 0$ and by proposition \ref{lim}, $\lim\limits_{(u,v)\to p}|H|=\lim\limits_{(u,v)\to p}\frac{|\HO|}{|\laO|}=\infty$. 

\item Since $\KO(\p)\neq 0$ and by proposition \ref{lim}, $\lim\limits_{(u,v)\to p}|K|=\lim\limits_{(u,v)\to p}\frac{|\KO|}{|\laO|}=\infty$.
		
\end{enumerate}
\end{proof}

\begin{lemma}\label{la}
	Let $\x:U \to \R^3$ be a proper wavefront, $W \subset U$ a compact set, $V \subset U$ an open set, $a:W\to\R$ a continuous function, $\Omegam_1$ and $\Omegam_2$ tangent moving bases of $\x$, then we have:
\begin{enumerate}
	\item If there exist a constant $C_1>0$ such that $|a|\leq C_1|\lambda_{\Omega_1}|$ on $W$ then there exist a constant $C_2>0$ such that $|a|\leq C_2|\lambda_{\Omega_2}|$ on $W$. 
	\item $\mathfrak{T}_{\Omega_1}(V)=\mathfrak{T}_{\Omega_2}(V)$ 
\end{enumerate}	
\end{lemma}
\begin{proof}
	\
	\begin{enumerate}
		\item Setting $\mathbf{A}=\I_{\Omega_1}^{-1}\Omegam_1^T\Omegam_2$ (change of basis matrix) and $\rho=det(\mathbf{A})$, we have $\Omegam_2=\Omegam_1 \mathbf{A}$, therefore $\Lambdam_{\Omega_1}=\mathbf{A}\Lambdam_{\Omega_2}$ and $\lambda_{\Omega_1}=\rho\lambda_{\Omega_2}$. Since $|a|\leq C_1|\lambda_{\Omega_1}|=C_1|\rho||\lambda_{\Omega_2}|$ and choosing $C_2$ as the maximum of $C_1|\rho|$ on $W$, we get the result. 
		\item Using the proof of item (1), $\lambda_{\Omega_1}=\rho\lambda_{\Omega_2}$ with $\rho\neq 0$, then we have the equality.   
	\end{enumerate}
\end{proof}

\begin{lemma}\label{lb}
	Let $\x_1:U \to \R^3$ be a proper wavefront with $U$ open connected, $W \subset U$ a compact set, $V \subset U$ an open set, $\mathbf{h}:Z\to U$ a diffeomorphism and $\Omegam$ a tmb of $\x$. Setting $\x_2:=\x_1\circ \mathbf{h}$ and choosing $\Omegam^h=(\Omegam \circ \mathbf{h})D\mathbf{h}$ as tmb of $\x_2$ we have:
	\begin{enumerate}
		\item There exist a constant $C_1>0$ such that $|e_1|\leq C_1|\lambda_{\Omega}|$, $|f_1|\leq C_1|\lambda_{\Omega}|$, $|g_1|\leq C_1|\lambda_{\Omega}|$ on $W$ if and only if there exist a constant $C_2>0$ such that $|e_2|\leq C_2|\lambda_{\Omega^h}|$, $|f_2|\leq C_2|\lambda_{\Omega^h}|$, $|g_2|\leq C_2|\lambda_{\Omega^h}|$ on $\mathbf{h}^{-1}(W)$. Where $e_1$, $f_1$, $g_1$ and $e_2$, $f_2$, $g_2$ are the coefficients of the second fundamental form of $\x_1$ and $\x_2$ respectively. 
		\item $e_1, f_1, g_1\in \mathfrak{T}_{\Omega}(V)$ if and only if $e_2, f_2, g_2\in\mathfrak{T}_{\Omega^h}(\mathbf{h}^{-1}(V))$. 
	\end{enumerate}	
\end{lemma}
\begin{proof} Let us denote by $\II_1$, $\II_2$ the matrices of the second fundamental forms of $\x_1$ and $\x_2$ respectively. If $det(D\mathbf{h})>0$, then $\n_2=\n_1\circ \mathbf{h}$ (in the case $det(D\mathbf{h})<0$, $\n_2=-\n_1\circ \mathbf{h}$ and it is analogous) therefore
$\II_2=-D\x_2^T D\n_2=-D\mathbf{h}^TD\x_1^TD\n_1D\mathbf{h}=D\mathbf{h}^T\II_1D\mathbf{h}$. This last equality expresses the coefficients $e_2$, $f_2$, $g_2$ as sum of multiples of the coefficients $e_1$, $f_1$, $g_1$ and vice versa. Since $\lambda_{\Omega^h}=\laO\circ\mathbf{h}$ (see proposition \ref{inv}) we get items (1) and (2) easily.    
\end{proof}
 
\begin{theorem}\label{K1}
Let $\x:U \to \R^3$ be a proper wavefront with just singularities of rank 1, $\Omegam$ a tmb of $\x$, $\p \in \Sigma(\x)$. Let $W\subset U$ be a compact neighborhood of $\p$ in which the relative principal curvature $k_{i\Omega}\neq0$ does not vanish. Let $k_j$ be the function that admit a $C^\infty$ extension to $W$ and $K$ the Gaussian curvature, then the following statements are equivalent: 
\begin{enumerate}
	\item $K$ is bounded on $W-\Sigma(\x)$.
	\item There exist a constant $C>0$ such that $|\KO|\leq C|\laO|$ on $W$.
	\item There exist a constant $C>0$ such that $|k_j|\leq C|\laO|$ on $W$.
	\item There exist a constant $C>0$ such that $|e|\leq C|\laO|$, $|f|\leq C|\laO|$ and $|g|\leq C|\laO|$ on $W$.  
\end{enumerate}	
\end{theorem}

\begin{proof}
	\
	\begin{itemize}
		\item $(1\Leftrightarrow 2)$ As $|K|=\frac{|\KO|}{|\laO|}$ on $W-\Sigma(\x)$ and by density of $W-\Sigma(\x)$ in $W$, it follows the equivalence.
		\item $(2\Leftrightarrow 3)$ Since $W$ is compact, $k_j$, $\KO$, $k_{i\Omega}$ are continuous on $W$ and    $k_j=\frac{\KO}{k_{i\Omega}}$, from this last equality follows the equivalence. 
		\item $(4\Rightarrow 1)$ we have that  $|eg-f^2|\leq 2C^2|\laO|^2$ on $W$, but by proposition \ref{D}, $EG-F^2=(\EO\GO-\FO^2)\laO^2$, then $|K|\leq \frac{2C^2}{\EO\GO-\FO^2}$ on $W-\Sigma(\x)$. Since $\frac{2C^2}{\EO\GO-\FO^2}$ is continuous on $U$ and $W$ is compact, $K$ is bounded on $W-\Sigma(\x)$.
		\item $(2\Rightarrow 4)$ if we prove $(4)$ locally on $W$, we can choose an open covering $B_k$ (open sets with the induced topology) of $W$ in which $(4)$ is satisfied in each compact $\bar{B_k}$ with constants $C_k$. Reducing this covering to a finite one, we have finite constants $C_{k_1},.., C_{k_n}$ and choosing $C$ as the maximum of these constants, $(4)$ is satisfied globally on $W$. 
		
		To prove this locally, first, for each $\q \in W$ let us take a tmb $D\y_l$ as in lemma \ref{lem1} on a neighborhood $V$ of $\q$ with $\y_l$ free of umbilical point (\ref{lem2}) on $V$. Shrinking $V$ if it is necessary, there exist a diffeomorphism $\mathbf{h}:V'\to V$, such that $\y_l\circ \mathbf{h}$ has derivatives as principal directions. By lemmas \ref{la} and \ref{lb}, we can assume that $\Omegam=D\y_{l}$ being $\y_l$ an embedding with derivatives as principal directions. Thus, by lemmas \ref{lem1} and \ref{lem2} $\I_{\Omega}$, $\II_{\Omega}$, $\boldsymbol{\alpha}_l$ and $\Lambdam=(\lambda_{ij})=\id-l\boldsymbol{\alpha}_l$ are diagonal matrices. If $rank(\Lambdam(\q))=1$, without loss of generality shrinking $V$ to a compact neighborhood, we can suppose that $\lambda_{22}(\q)=1-l\alpha_{l22}(\q)=0$, with $\lambda_{11}\neq0$ and $-\frac{\go}{\GO}=\alpha_{l22}\neq0$ on $V$. By proposition \ref{D}, $f=0$, $e=\lambda_{11}\eo$, $g=\lambda_{22}\go=\laO\frac{\go}{\lambda_{11}}$ and by hypothesis $(2)$ $|\frac{\go\eo}{\EO\GO}|\leq C|\laO|$ on $V'=V\cap W$, thus $|e|\leq C|\lambda_{11}||\frac{\EO\GO}{\go}||\laO|$. If we choose $C'$ as the biggest maximum of the functions $C|\lambda_{11}||\frac{\EO\GO}{\go}|$ and $|\frac{\go}{\lambda_{11}}|$ on $V'$, we get that $|e|\leq C'|\laO|$, $|f|\leq C'|\laO|$ and $|g|\leq C'|\laO|$ on $V'$.
		
		On the other hand, if $rank(\Lambdam(\q))=2$, shrinking $V$ to a compact neighborhood, we can suppose that $\lambda_{11}\neq0$ and $\lambda_{22}\neq0$ on $V$. Thus, $f=0$, $e=\lambda_{11}\eo=\laO\frac{\eo}{\lambda_{22}}$, $g=\lambda_{22}\go=\laO\frac{\go}{\lambda_{11}}$, then choosing $C'$ as the biggest maximum of the functions $|\frac{\eo}{\lambda_{22}}|$ and $|\frac{\go}{\lambda_{11}}|$ on $V'=V\cap W$, we have $|e|\leq C'|\laO|$, $|f|\leq C'|\laO|$ and $|g|\leq C'|\laO|$ on $V'$.            
	\end{itemize}
\end{proof}

\begin{ex}\label{ex1}
The wavefront $\x(u,v)=(u, 2v^3+u^2v,3v^4+u^2v^2)$ (cuspidal lips) with normal vector $\n=(2uv^2,-2v,1)(4u^2v^4+4v^2+1)^{-\frac{1}{2}}$ has a singularity at $(0,0)$, Gaussian curvature $K=-\frac{4v^2}{(4u^2v^4+4v^2+1)^2(u^2+6v^2)}$ with $|K|\leq 1$. Observe that $K$ does not converge when $(u,v)\longrightarrow (0,0)$, then it is not extendable.
\begin{figure}[h]
\begin{center}
\includegraphics[scale=0.50]{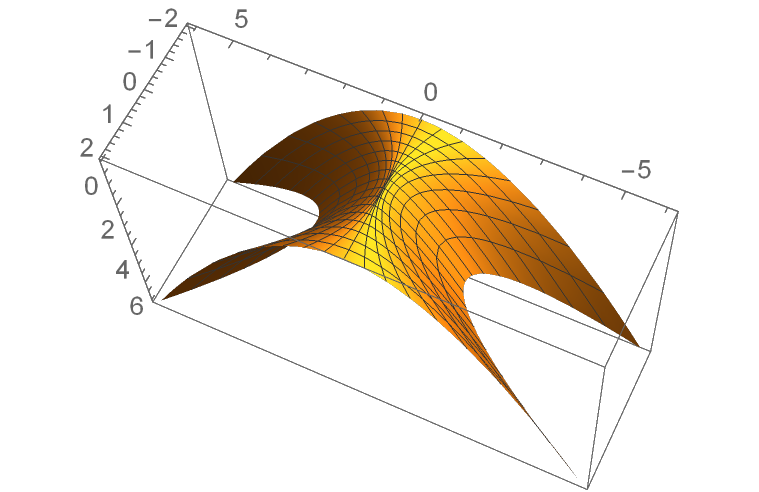}\qquad\qquad \includegraphics[scale=0.50]{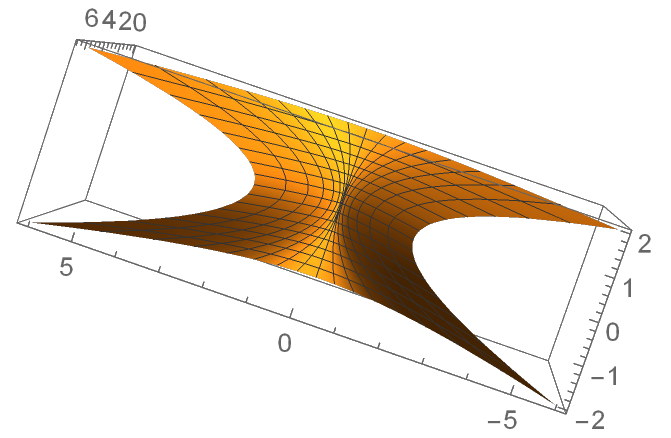}
\end{center}
\caption{A wavefront with degenerate singularity of rank $1$ at the origin and Gaussian curvature bounded but non-extendable.}\label{bound}
\end{figure}
\end{ex}

\begin{theorem}\label{K2}
Let $\x:U \to \R^3$ be a proper wavefront with just singularities of rank 1, $\Omegam$ a tmb of $\x$, $\p \in \Sigma(\x)$ with $rank(D\x(\p))=1$. Let $V\subset U$ be a open neighborhood of $\p$ in which the relative principal curvature $k_{i\Omega}\neq0$ does not vanish. Let $k_j$ be the function that admit a $C^\infty$ extension to $V$ and $K$ the Gaussian curvature, then the following statements are equivalent: 
\begin{enumerate}
	\item The Gaussian curvature $K$ admits a $C^\infty$ extension to $V$.
	\item $\KO \in \mathfrak{T}_\Omega(V)$.
	\item $k_j \in \mathfrak{T}_\Omega(V)$.
	\item $e, f, g \in \mathfrak{T}_\Omega(V)$.
\end{enumerate}	
where $\mathfrak{T}_\Omega(V)$ is the principal ideal generated by $\laO$ in the ring $C^\infty(V,\R)$
\end{theorem}
\begin{proof}
	\
	\begin{itemize}
		\item $(1\Leftrightarrow 2)$ As $\KO=K\laO$ on $V-\Sigma(\x)$ and by density of $V-\Sigma(\x)$ in $V$, it follows the equivalence.
		\item $(2\Leftrightarrow 3)$ Since    $k_j=\frac{\KO}{k_{i\Omega}}$, from this last equality follows the equivalence.
		\item $(4\Rightarrow 1)$ we have that  $eg-f^2=\phi\laO^2$ with $\phi \in C^\infty(V,\R)$, but by proposition \ref{D}, $EG-F^2=(\EO\GO-\FO^2)\laO^2$, then $K= \frac{\phi}{\EO\GO-\FO^2}$ on $V-\Sigma(\x)$. Since $\frac{\phi}{\EO\GO-\FO^2}$ is smooth on $V$, $K$ has a $C^\infty$ extension to $V$.
		\item $(2\Rightarrow 4)$ if we prove $(4)$ locally on $V$, we can choose an open covering locally finite $B_k\subset V$ (open balls) of $V$, $k\in \mathbb{N}$ with a partition of the unity $\psi_k$ subordinated to this open cover in which $e, f, g \in \mathfrak{T}_\Omega(B_k)$ for every $B_k$. For each $k\in \mathbb{N}$ there exist $f_{1k}, f_{2k}, f_{3k}\in C^\infty(V,\R)$ such that $e=f_{1k}\laO$, $f=f_{2k}\laO$, $g=f_{3k}\laO$ on $B_k$. Since that the supports of $f_{sk}\psi_k$ form families locally finite for $s= 1 ,2, 3$, we have that $f_s:=\sum_k f_{sk}\psi_k \in C^\infty(V,\R)$ for $s= 1 ,2, 3$, therefore $e=f_{1}\laO$, $f=f_{2}\laO$, $g=f_{3}\laO$ on $V$.  
		
		To prove this locally, first, for each $\q \in V$ let us take a tmb $D\y_l$ as in lemma \ref{lem1} on a neighborhood $Z\subset V$ of $\q$ with $\y_l$ free of umbilical point on $Z$. Shrinking $Z$ if it is necessary, there exist a diffeomorphism $\mathbf{h}:Z'\to Z$, such that $\y_l\circ \mathbf{h}$ has derivatives as principal directions. By lemmas \ref{la} and \ref{lb}, we can assume that $\Omegam=D\y_{l}$ being $\y_l$ an embedding with derivatives as principal directions on $Z$. Thus, by lemmas \ref{lem1} and \ref{lem2} $\I_{\Omega}$, $\II_{\Omega}$, $\boldsymbol{\alpha}_l$ and $\Lambdam=(\lambda_{ij})=\id-l\boldsymbol{\alpha}_l$ are diagonal matrices. If $rank(\Lambdam(\q))=1$, without loss of generality shrinking $Z$ to a open neighborhood $V'$, we can suppose that $\lambda_{22}(\q)=1-l\alpha_{l22}(\q)=0$, with $\lambda_{11}\neq0$ and $-\frac{\go}{\GO}=\alpha_{l22}\neq0$ on $V'$. By proposition \ref{D}, $f=0$, $e=\lambda_{11}\eo$, $g=\lambda_{22}\go=\laO\frac{\go}{\lambda_{11}}$ and by hypothesis $(2)$ $\frac{\go\eo}{\EO\GO}= \phi\laO$ for some $\phi \in C^\infty(V,\R)$, then $e= \phi\lambda_{11}\frac{\EO\GO}{\go}\laO$. Thus, we get that $e\in \mathfrak{T}_\Omega(V')$, $f\in \mathfrak{T}_\Omega(V')$ and $g\in \mathfrak{T}_\Omega(V')$.
		
		On the other hand, if $rank(\Lambdam(\q))=2$, shrinking $Z$ to a open neighborhood $V'$, we can suppose that $\lambda_{11}\neq0$ and $\lambda_{22}\neq0$ on $V'$. Thus, $f=0$, $e=\lambda_{11}\eo=\laO\frac{\eo}{\lambda_{22}}$, $g=\lambda_{22}\go=\laO\frac{\go}{\lambda_{11}}$, then we have $e\in \mathfrak{T}_\Omega(V')$, $f\in \mathfrak{T}_\Omega(V')$ and $g\in \mathfrak{T}_\Omega(V')$.
	\end{itemize}
\end{proof}

Let $\mathcal{E}_{2,p}$ denote the local ring of smooth function germs at the point $\p$. By $[g]_p$ we denote the function germ at $\p$ with $g \in C^\infty(U,\R)$ a representative of the germ defined on a neighborhood $U$ of $\p$. We define $J_p$ as the ideal in $\mathcal{E}_{2,p}$ generated by $[\laO]_p$,  $\hat{J_p}:=\{[g]_p\in \mathcal{E}^1_{2,p}: \textnormal{there exist $C>0$ such that  $|g|\leq C|\laO|$ on some neighborhood of $\p$}\}$ and $J_{\Sigma_p }:=\{[g]_p\in \mathcal{E}^1_{2,p}: \textnormal{for some neighborhood $U$ of $\p$, $g$ vanish on $U\cap \laO^{-1}(0)$}\}$. These ideals satisfy $ J_p\subset \hat{J_p}\subset J_{\Sigma_p}$, their definitions do not depend on the chosen tmb $\Omegam$ and when $\p$ is a non-degenerate singularity these three ideals are equal. To see that, we can assume that $\p=\mathbf{0}$ and making a change of coordinates, we can assume that $\laO$ is equal to $u$ or $v$, then applying the Hadamard lemma \cite[ Ch IV, page 100]{Gib}, we obtain the result. From this and theorems \ref{K1} and \ref{K2} we have the following corollary.

\begin{corollary}

	Let $\x:U \to \R^3$ be a proper wavefront, $\Omegam$ a tmb of $\x$, $\p \in \Sigma(\x)$ a non-degenerate singularity with $rank(D\x(\p))=1$. Let $k_j$ be the function that admit a local $C^\infty$ extension at $\p$ and $K$ the Gaussian curvature, then the following statements are equivalent: 
	\begin{enumerate}
		\item The Gaussian curvature $K$ admits a local $C^\infty$ extension at $\p$.
		\item The Gaussian curvature $K$ is locally bounded on some neighborhood of $\p$.
		\item $[\KO]_p \in J_{\Sigma p}$.
		\item $[k_j]_p \in J_{\Sigma p}$.
		\item $[e]_p, [f]_p, [g]_p \in J_{\Sigma p}$.
	\end{enumerate}	
Where $J_{\Sigma p}=\hat{J_p}=J_p$ in this case.	

\end{corollary}
The equivalence between $(1)$, $(2)$ and $(5)$ were obtained by K. Saji, M. Umehara, and K. Yamada in \cite{front}.

\begin{ex}\label{ex2}
	The wavefront $(u,\sin(ku)\frac{v^{k+1}}{k+1},\sin(ku)\frac{v^{k+2}}{k+2})$, $k\in \mathbb{N}$, has as a tangent moving basis:
	$$
	\Omegam=\begin{pmatrix}
	1&0\\
	\cos(ku)k\frac{v^{k+1}}{k+1}&1\\
	\cos(ku)k\frac{v^{k+2}}{k+2}&v
	\end{pmatrix}, \ \Lambdam=\begin{pmatrix}
	1&0\\
	0&\sin(ku)v^k
	\end{pmatrix}$$  $$\II_\Omega=\begin{pmatrix}
	\sin(ku)k^2v^{k+2}(\frac{1}{k+1}-\frac{1}{k+2})&0\\
	0&1
	\end{pmatrix}\frac{1}{\sqrt{\epsilon}}$$
	where $\epsilon=1+v^2+\cos^2(ku)k^2v^{2k+4}(\frac{1}{k+1}-\frac{1}{k+2})^2$. Then, $\laO=\sin(ku)v^k$, $(\eo\go-\fuo\fdo) \in \mathfrak{T}_\Omega$ therefore $\KO \in \mathfrak{T}_\Omega$ and by the proposition above the Gaussian curvature $K$ admit a $C^\infty$ extension to $\R^2$. Observe that, since $\II=\Lambdam\II_\Omega$ (see proposition \ref{D}), we have $e, f, g \in \mathfrak{T}_\Omega$.   
\end{ex}

\begin{remark}
	The boundedness and extendibility of Gaussian curvature are conserved under changes of coordinates in the domain, however they are not conserved making changes at co-domain. The wavefront $\x(u,v)=(u, 2v^3+u^2v,3v^4+u^2v^2+u^2)$ with Gaussian curvature unbounded can be obtained from example \ref{ex1} whose Gaussian curvatures is bounded, applying at co-domain the diffeomorphism $F(X,Y,Z)=(X,Y,Z+X^2)$. The same situation occurs with the wavefronts from example \ref{ex2} and $(u,\sin(ku)\frac{v^{k+1}}{k+1},\sin(ku)\frac{v^{k+2}}{k+2}+u^2)$ which have extendable and  non-extendable Gaussian curvatures respectively. 
\end{remark}

\begin{definition}
	Let $\x:U \to \R^3$ be a wavefront, $\Omegam$ a tmb of $\x$ and $\p \in \Sigma(\x)$. We say that $\x$ is {\it parallelly smoothable at $\p$} if there exist $\epsilon>0$ and an open neighborhood $V$ of $\p$ such that $rank(D(\x+l\n)(\q))=2$ for every $(\q,l) \in V\times(0,\epsilon)$ or every $(\q,l) \in V\times(-\epsilon,0)$.
\end{definition}

\begin{theorem}\label{pa1}
Let $\x:U \to \R^3$ be a wavefront, $\Omegam$ a tmb of $\x$, $\p \in \Sigma(\x)$ with $rank(D\x(\p))=1$, $\KO(\p)\neq 0$ and $\HO(\p)<0$ (resp. $\HO(\p)>0$) then the following statements are equivalents:
\begin{enumerate}
\item $\x$ is parallelly smoothable at $\p$.
\item $\lim\limits_{(u,v)\to p}H=\pm\infty$.
\item $\lim\limits_{(u,v)\to p}k_1=\pm\infty$ (resp. $k_2$).
\item $\lim\limits_{(u,v)\to p}K=\pm\infty$.
\item There exist an open neighborhood $V$ of $\p$ in which $\laO$ does not change sign.
\item There exist an open neighborhood $V$ of $\p$ in which $k_{2\Omega}$ (resp. $k_{1\Omega}$) does not change sign.  
\end{enumerate}
	
\end{theorem}

\begin{proof}
	\
	\begin{itemize}
		\item $(1\Rightarrow5)$ if $\x$ is parallelly smoothable at $\p$ if and only if there exist $\epsilon>0$ and an open neighborhood $V$ of $\p$ such that $\y_t=\x+t\n|_V$ is an immersion for every $t\in (0,\epsilon)$ (or every $t\in(-\epsilon,0)$) if and only if $det(\Lambdam^T+t\boldsymbol{\mu}^T)(\q)=\laO(\q)-2t\HO(\q)+t^2\KO(\q)\neq0$ for every $t\in (0,\epsilon)$ (or every $t\in(-\epsilon,0)$) and $\q \in V$. Shrinking $V$ we can suppose this is connected, then we have that $\laO(\q)-2t\HO(\q)+t^2\KO(\q)>0$ (or $<0$) for every $(\q,t)\in V\times(0,\epsilon)$, thus taking the limit in both sides of this inequality when $t$ tends to $0$, we get that $\laO\geq0$ on $V$.
		\item $(5\Leftrightarrow2)$ as $H=\frac{\HO}{\laO}$ on $U-\Sigma(\x)$ and $\HO(\p)\neq0$, follows the equivalence.
		\item $(5\Leftrightarrow4)$ as $K=\frac{\KO}{\laO}$ on $U-\Sigma(\x)$ and $\KO(\p)\neq0$, follows the equivalence.
		\item $(2\Leftrightarrow3)$ by proposition \ref{prin1}, there exist a neighborhood $V$ of $\p$ such that $k_2$ has a $C^\infty$ extension and since $k_1=2H-k_2$ follows the equivalence.
		\item $(3\Leftrightarrow6)$ there exist a neighborhood $W$ of $\p$ such that $\HO<0$ and $\KO\neq0$ on $W$, then $k_1\neq0$ on $W-\Sigma(\x)$. Since $\frac{\KO}{k_1}=k_{2\Omega}$ on $W-\Sigma(\x)$ and $k_{2\Omega}=0$ on $\Sigma(\x)$(by proposition \ref{prin1}), follows the equivalence.
		\item $(6\Rightarrow1)$ there exist a neighborhood $W$ of $\p$ such that $\HO<0$ and $\KO\neq0$ on $W$. By proposition \ref{prin1} $k_{1\Omega}(\p)\neq0$, thus $\frac{k_{1\Omega}(\p)}{\KO(\p)}\neq0$, then using item $(6)$ there exist $\epsilon>0$ and an open connected $V$ of $\p$ such that $\frac{k_{2\Omega}}{\KO}$ does not change sign and $|\frac{k_{1\Omega}}{\KO}|>\epsilon$ on $V$. Thus, if $\frac{k_{2\Omega}}{\KO}\geq0$ (resp. $\frac{k_{2\Omega}}{\KO}\leq0$) then $\laO(\q)-2t\HO(\q)+t^2\KO(\q)\neq0$ for every $(\q, t)\in V\times (-\epsilon,0)$ (resp. $V\times (0,\epsilon)$) because $\laO(\q)-2t\HO(\q)+t^2\KO(\q)=0$ if and only if $t=\frac{k_{1\Omega}}{\KO}(\q)$ or $t=\frac{k_{2\Omega}}{\KO}(\q)$.  
	\end{itemize}
	
\end{proof}

\begin{theorem}\label{pa2}
	Let $\x:U \to \R^3$ be a wavefront, $\Omegam$ a tmb of $\x$, $\p \in \Sigma(\x)$ with $rank(D\x(\p))=1$, $\KO(\p)=0$ and $\HO(\p)<0$ (resp. $\HO(\p)>0$) then the following statements are equivalents:
	\begin{enumerate}
		\item $\x$ is parallelly smoothable at $\p$.
		\item $\lim\limits_{(u,v)\to p}H=\pm\infty$.
		\item $\lim\limits_{(u,v)\to p}k_1=\pm\infty$(resp. $k_2$).
		\item There exist an open neighborhood $V$ of $\p$, in which $\laO$ does not change sign.  
	\end{enumerate}
	
\end{theorem}

\begin{proof} The proof of $(1\Rightarrow4)$, $(2\Leftrightarrow3)$, $(2\Leftrightarrow4)$ is equal to corresponding ones in theorem \ref{pa1} because this does not use the hypothesis of $\KO(\p)\neq0$. To prove $(3\Rightarrow1)$, let us define $A:=\{\q\in U: \KO(\q)\neq0\}$. If $\p \notin \bar{A}$, there exist a neighborhood $W$ of $\p$ in which $\KO\equiv0$, thus $\laO(\q)-2t\HO(\q)+t^2\KO(\q)=\laO(\q)-2t\HO(\q)$ on $W$ and since $\HO(\p)\neq0$, using that $\laO$ does not change sign on a neighborhood of $\p$, shrinking $W$ we have that $\frac{\laO}{\HO}$ does not change sign, then if $\frac{\laO}{\HO}\geq0$ (resp. $\frac{\laO}{\HO}\leq0$) we have that $\laO(\q)-2t\HO(\q)\neq0$ for every $(\q, t)\in W\times (-\epsilon,0)$ (resp. $W\times (0,\epsilon)$) with an arbitrary $\epsilon>0$ and from this follows the result. If $\p \in \bar{A}$, by hypothesis there exist an open neighborhood $W$ of $\p$ such that $k_1>0$ and $H>0$ (resp. or $<0$) on $W-\Sigma(\x)$ and as $k_{1\Omega}(\p)\neq0$, we have that $\lim\limits_{(u,v)\to p}|\frac{k_{1\Omega}}{\KO}|\big|_{A}=\infty$. Let $\epsilon>0$ be given, there exist a open ball $B$ such that $|\frac{k_{1\Omega}}{\KO}|>\epsilon$ on $B\cap A$, $k_1>0$ (resp. $<0$), $H>0$ (resp. $<0$) and $\frac{\laO}{\HO}=\frac{1}{H}>0$ (resp. $<0$) on $B-\Sigma(\x)$. Since $\frac{k_{2\Omega}}{\KO}=\frac{1}{k_1}$ on $(B\cap A)-\Sigma(\x)$, we have that $\frac{k_{2\Omega}}{\KO}\geq0$ (resp. $\leq0$) on $B\cap A$. Now, if $(\q, t)\in B\times (-\epsilon,0)$ (resp. $B\times (0,\epsilon)$) and $\KO(\q)=0$ then $\laO(\q)-2t\HO(\q)+t^2\KO(\q)=\laO(\q)-2t\HO(\q)\neq0$ because  $\frac{\laO}{\HO}\geq0$ (resp. $\leq0$) on $B$. The another option is that $\KO(\q)\neq0$, then $\laO(\q)-2t\HO(\q)+t^2\KO(\q)\neq0$ because this is $0$ if and only if $t=\frac{k_{1\Omega}}{\KO}(\q)$ or $t=\frac{k_{2\Omega}}{\KO}(\q)$ which is impossible since that $|\frac{k_{1\Omega}}{\KO}(\q)|>\epsilon$ and $\frac{k_{2\Omega}}{\KO}(\q)\geq0$ (resp. $\leq0$). It follows the result.          
\end{proof}

\begin{corollary}\label{app}
	Let $\x:U \to \R^3$ be a wavefront, $\Omegam$ a tmb of $\x$, $\p \in \Sigma(\x)$ with $rank(D\x(\p))=1$, we have that $\x$ is parallelly smoothable at $\p$ if and only if $\laO$ does not change sign on a neighborhood of $\p$.
\end{corollary}

\begin{corollary}
	Let $\x:U \to \R^3$ be a wavefront, $\Omegam$ a tmb of $\x$, $\p \in \Sigma(\x)$ with $rank(D\x(\p))=1$, if $\x$ is parallelly smoothable at $\p$, then $\p$ is a degenerate singularity.
\end{corollary}
\begin{proof}
	If we suppose that $\p=(p_1,p_2)$ is a non-degenerate singularity, then $\lambda_{\Omega u}(\p)\neq 0$ or $\lambda_{\Omega v}(\p)\neq 0$ and therefore $\laO(u,p_2)$ or $\laO(p_1,v)$ is strictly monotone as function of one variable on every sufficient small neigborhood of $\p$, which is contradictory, because $\laO(\p)=0$ and this does not change sign by corollary \ref{app}.   
\end{proof}

\begin{corollary}
	Let $\x:U \to \R^3$ be a wavefront, $\Omegam$ a tmb of $\x$, $\p \in \Sigma(\x)$ with $rank(D\x(\p))=1$, if $\p$ is an isolated singularity then  $\x$ is parallelly smoothable at $\p$.
\end{corollary}
\begin{proof}
	if $\p$ is an isolated singularity, then there exist an open connected neighborhood $V$ of $\p$, such that $\laO\neq 0$ on $V-\{\p\}$ and since that $V-\{\p\}$ is arc-connected, $\laO$ does not change sign on $V$. By corollary \ref{app}, it follows the result.  
\end{proof}

\begin{ex}

	The wavefront $(u, 2v^3+u^2v,3v^4+u^2v^2)$ (cuspidal lips) has a isolated singularity at $(0,0)$, then it is parallelly smoothable at $(0,0)$. On the other hand, $(u, 2v^3-u^2v,3v^4-u^2v^2)$ (cuspidal beaks) is not parallelly smoothable at $(0,0)$, because taking as a tangent moving basis:
	$$
	\Omegam=\begin{pmatrix}
	1&0\\
	-2uv&1\\
	-2uv^2&2v
	\end{pmatrix}, \ \Lambdam=\begin{pmatrix}
	1&0\\
	0&6v^2-u^2
	\end{pmatrix}$$
we get $\laO=6v^2-u^2$, which changes of sign on every neighborhood of $(0,0)$. By the same argument, $\x(u,v)=(u, v^2, v^3)$ (cuspidal edge) and $\x(u,v)=(3u^4 + u^2v, 4u^3 + 2uv, v)$ (swallowtail) are not parallelly smoothable at $(0,0)$, because can be chosen tmb's $\Omegam$ in which $\laO$ is $v$ and $12u^2+2v$ respectively.
\end{ex}

\begin{theorem}\label{ka}
	Let $\x:U \to \R^3$ be a wavefront, $\Omegam$ a tmb of $\x$ and $\p \in \Sigma(\x)$ with $rank(D\x(\p))=1$. One of the principal curvatures $\kappa_{-}, \kappa_{+}$ has a $C^\infty$-extension to an open neighborhood of $\p$ if and only if $\x$ is parallelly smoothable at $\p$.
\end{theorem}
\begin{proof}
	if $\x$ is parallelly smoothable at $\p$, then by corollary \ref{app} $\laO\geq0$ (or $\laO\leq0$, this case is analogous) on an open neighborhood $V$ of $\p$, thus $k_1=\kappa_{-}$, $k_2=\kappa_{+}$ on $V-\Sigma(\x)$ and by proposition \ref{prin1} one of these function has a $C^\infty$-extension to an open neighborhood of $\p$. Conversely, without loss of generality let us suppose that $\kappa_{-}$ has a $C^\infty$-extension to an open neighborhood $W$ of $\p$, then $\laO$ does not change sign on some neighborhood of $\p$, otherwise there are sequences $\mathbf{a}_n\longrightarrow\p$, $\mathbf{b}_n\longrightarrow\p$ such that $\laO(\mathbf{a}_n)>0$ and $\laO(\mathbf{b}_n)<0$ for every $n\in \mathbb{N}$. Thus, $\lim\limits_{n\to \infty}|k_1(\mathbf{a}_n)|=\lim\limits_{n\to \infty}|\kappa_{-}(\mathbf{a}_n)|=|\kappa_{-}(\p)|=\lim\limits_{n\to \infty}|\kappa_{-}(\mathbf{b}_n)|=\lim\limits_{n\to \infty}|k_2(\mathbf{b}_n)|$ which is contradictory, because by proposition \ref{prin1} one of the limits $\lim\limits_{n\to \infty}|k_1(\mathbf{a}_n)|$, $\lim\limits_{n\to \infty}|k_2(\mathbf{b}_n)|$ is $\infty$.         
\end{proof}

\begin{remark}
Observe that, if $\kappa_{-}$ (or $\kappa_{+}$) have a $C^\infty$-extension locally at $\p$, then the another one converge to $\pm\infty$ at $\p$ by items (3)s of theorems \ref{pa1} and \ref{pa2}. In fact, reasoning similarly as in theorem \ref{ka}, we have that one of the principal curvatures $\kappa_{-}, \kappa_{+}$ converge to $\pm\infty$ at $\p$ if and only if $\x$ is parallelly smoothable at $\p$. Also can be proved similarly that $\kappa_{-}$ (resp. $\kappa_{+}$) is bounded locally at $\p$ if and only if $\kappa_{-}$ (resp. $\kappa_{+}$) have a $C^\infty$-extension locally at $\p$.   
\end{remark}

\section{Singularities of rank 0}

\begin{proposition}\label{propc2}
	Let $\x:U \to \R^3$ be a proper wavefront, $\Omegam$ a tmb of $\x$, $\p \in \Sigma(\x)$ with $rank(D\x(\p))=0$, then:
	\begin{enumerate}
		\item $(k_{1\Omega},k_{2\Omega})(\p)=(0,0)$.
		\item $\lim\limits_{(u,v)\to p}|K|=\infty$.
		\item $\frac{1}{k_1}$ and $\frac{1}{k_2}$ have  continuous extensions on a neighborhood $V$ of $\p$, which are of class $C^{\infty}$ except possibly at umbilical points and singularities of rank 0 of $\x$.  
		\item $\lim\limits_{(u,v)\to p}|k_1|=\infty$ and $\lim\limits_{(u,v)\to p}|k_2|=\infty$.
		\item $\frac{1}{k_1}+\frac{1}{k_2}$ has a $C^\infty$-extension on a neighborhood $V$ of $\p$.
	\end{enumerate}
\end{proposition}

\begin{proof}
	\
	\begin{enumerate}
		\item By proposition \ref{wft}, $\HO(\p)=0$, then  $k_{1\Omega}(\p)=\HO(\p)-\sqrt{\HO^2(\p)-\laO(\p)\KO(\p)}=0$. Similarly $k_{2\Omega}(\p)=0$.
		\item By proposition \ref{wft}, $\KO(\p)\neq0$ and therefore $\lim\limits_{(u,v)\to p}|K|=\lim\limits_{(u,v)\to p}|\frac{\KO}{\laO}|=\infty$
		\item There exist a neighborhood $V$ of $\p$ such that $\KO(\p)\neq0$ and hence $k_1\neq0, k_2\neq0$ on $V-\Sigma(\x)$, then $\frac{1}{k_1}=\frac{k_2}{k_1k_2}=\frac{k_{2\Omega}}{\KO}$ and similarly $\frac{1}{k_2}=\frac{k_{1\Omega}}{\KO}$ which are well defined on $V$. Notice that $k_{1\Omega}$ and $k_{2\Omega}$ may not be differentiable at umbilical points and singularities of rank 0 of $\x$.  
		\item By item (2), there exist a neighborhood $V$ of $\p$ such that $k_1\neq0, k_2\neq0$ on $V-\Sigma(\x)$, therefore $k_{1\Omega}\neq0, k_{2\Omega}\neq0$ as well. Then, $k_1=\frac{\KO}{k_{2\Omega}}$ and $k_2=\frac{\KO}{k_{1\Omega}}$ on $V-\Sigma(\x)$ and using that $\KO(\p)\neq0$ and item (1) we get (3).
		\item Since $\KO(\p)\neq 0$, there exist an open neighborhood $V$ of $\p$ such that $\KO$ does not vanish on $V$, then $\frac{1}{k_1}+\frac{1}{k_2}=\frac{2H}{K}=\frac{2\HO}{\KO}$ on $V-\Sigma(\x)$ and as $\frac{2\HO}{\KO}$ is well defined on $V$, this is a $C^{\infty}$-extension. 
	\end{enumerate}
\end{proof}

\begin{theorem}\label{p2}
	Let $\x:U \to \R^3$ be a proper wavefront, $\Omegam$ a tmb of $\x$, $\p \in \Sigma(\x)$ with $rank(D\x(\p))=0$, then the following statements are equivalent:
	\begin{enumerate}
		\item $\x$ is parallelly smoothable at $\p$.
		\item $\laO\KO\geq 0$ and $\HO$ does not change sign on a neighborhood $V$ of $\p$.
		\item $\lim\limits_{(u,v)\to p}K=\infty$ and $\lim\limits_{(u,v)\to p}H=\pm\infty$.
		\item $\lim\limits_{(u,v)\to p}k_1=\lim\limits_{(u,v)\to p}k_2=\infty$ or $\lim\limits_{(u,v)\to p}k_1=\lim\limits_{(u,v)\to p}k_2=-\infty$.
	\end{enumerate}
\end{theorem}

\begin{proof}
	\
	\begin{itemize}
		\item $(1\Rightarrow2)$ If $\x$ is parallelly smoothable at $\p$, then there exist $\epsilon>0$ and an open connected neighborhood $V$ of $\p$ such that $\laO(\q)-2t\HO(\q)+t^2\KO(\q)\neq0$ for every $(\q,t)\in V\times(0,\epsilon)$ (or $V\times(-\epsilon,0)$, this case is analogous). As $\KO\neq0$ and does not change sign on $V$ (shrinking $V$ if it is necessary), $k_{1\Omega}=k_{2\Omega}=0$ and since $\laO(\q)-2t\HO(\q)+t^2\KO(\q)\neq0$ if and only if $t=\frac{k_{1\Omega}}{\KO}(\q)$ or $t=\frac{k_{2\Omega}}{\KO}(\q)$, we have that $\frac{k_{1\Omega}}{\KO}\leq0$ and $\frac{k_{2\Omega}}{\KO}\leq0$ on $V$. Then, $k_{1\Omega}k_{2\Omega}\geq0$ on $V$, but $-k_{1\Omega}$ and $-k_{2\Omega}$ are the eigenvalues of $\boldsymbol{\alpha}_\Omega=\boldsymbol{\mu}_{\Omega}adj(\Lambdam)$, then $\laO\KO=k_{1\Omega}k_{2\Omega}\geq0$ on $V$. Observe that $k_{1\Omega}$  and $k_{2\Omega}$ do not change sign on $V$, then $\HO$ neither.   
		\item $(2\Rightarrow3)$ Since $\laO\KO=\laO^2K$ on $U-\Sigma(\x)$ and using that $\lim\limits_{(u,v)\to p}|K|=\infty$ we get that $\lim\limits_{(u,v)\to p}K=\infty$. On the other hand, $H^2\geq K$, then $\lim\limits_{(u,v)\to p}|H|=\infty$. As $\HO$ and $\laO$ do not change sign on a neighborhood of $\p$, $H=\frac{\HO}{\laO}$ neither and we get the result.   
		\item $(3\Rightarrow4)$  As $K$ is positive near to $\p$ then $\laO\KO=\laO^2K\geq0$ and $\KO\neq0$ on a neighborhood $Z$ of $\p$. Shrinking $Z$, $\laO$ does not change sign and $H$ neither on $Z-\Sigma(\x)$. Without loss of generality, let us suppose $\laO\geq0$ on $Z$, then $k_1=H-\sqrt{H^2-K}$ and $k_2=H+\sqrt{H^2-K}$ and since that $H$ does not change sign on $Z-\Sigma(\x)$, one of the function $k_1, k_2$ neither. By this last and since that $K>0$, we have that $k_1>0, k_2>0$ or $k_1<0, k_2<0$ on $Z-\Sigma(\x)$, then using item (3) of proposition \ref{propc2} we get the result.      
		\item $(4\Rightarrow1)$ There exist a neighborhood $V$ of $\p$ such that $k_1>0, k_2>0$ (or $k_1<0, k_2<0$, this case is analogous) on $V-\Sigma(\x)$ and $\KO\neq0$ on $V$, then $\frac{k_{1\Omega}}{\KO}=\frac{1}{k_2}>0$ and $\frac{k_{2\Omega}}{\KO}=\frac{1}{k_1}>0$ on $V-\Sigma(\x)$, thus by density of $V-\Sigma(\x)$ $\frac{k_{1\Omega}}{\KO}, \frac{k_{2\Omega}}{\KO}\geq0$ on $V$. Choose $\epsilon>0$ arbitrary and we have that $\laO(\q)-2t\HO(\q)+t^2\KO(\q)\neq0$ for every $(\q,t)\in V\times(-\epsilon,0)$. It follows (1).   
		      
	\end{itemize}
	
\end{proof}

\begin{corollary}\label{app2}
	Let $\x:U \to \R^3$ be a wavefront, $\Omegam$ a tmb of $\x$, $\p \in \Sigma(\x)$ with $rank(D\x(\p))=0$. If there exist a neighborhood $V$ of $\p$ in which the only singularity of rank $0$ is $\p$, then $\x$ is parallelly smoothable at $\p$ if and only if $\laO\KO\geq 0$ on a neighborhood $W$ of $\p$.  
\end{corollary}
\begin{proof}
	If $\laO\KO\geq 0$ on $W$, shrinking if it is necessary we can suppose that $K\neq 0$ on $W-\Sigma(\x)$ and since $\laO\KO=\laO^2K$, then $H^2>K>0$ on $W-\Sigma(\x)$. As $\HO\neq0$ on singularities of rank 1, then $\HO$ has a isolated zero on $W\cap V$ and therefore $\HO$ does not change sign. Applying the last theorem we get the result.  
\end{proof}

Observe that, if we have a wavefront $\x:U \to \R^3$, $\Omegam$ a tmb of $\x$, $\p \in \Sigma(\x)$ with $rank(D\x(\p))=0$ and this is parallelly smoothable at $\p$, since that $\KO(\p)\neq0$ and $\laO\KO\geq 0$ on a neighborhood of $\p$, then $\laO$ does not change sign on a neighborhood of $\p$. However, this condition is not sufficient as it happened in the case of singularities of rank 1. The next example shows this.

\begin{ex}
	The wavefront $\x:=(u^k,\pm v^k,\frac{k}{k+1}u^{k+1}\pm\frac{k}{k+1}v^{k+1})$, with $k\in \mathbb{N}$, $k\geq2$ has as tmb:
	$$
	\Omegam=\begin{pmatrix}
	1&0\\
	0&1\\
	u&v
	\end{pmatrix}, \ \Lambdam=\begin{pmatrix}
	ku^{k-1}&0\\
	0&\pm kv^{k-1}
	\end{pmatrix},$$
then $\laO=\pm k^2u^{k-1}v^{k-1}$, $\KO=\frac{1}{(1+u^2+v^2)^2}$. By corollary \ref{app2}, $\x$ is parallelly smoothable at $(0,0)$ when we choose $k$ odd and the sign $+$. If $k$ is even or the sign is $-$, this is not parallelly smoothable at $(0,0)$, even when $k$ is odd with sign $-$ in the expression, in which $\laO$ does not change sign.  	
\end{ex}

\begin{corollary}
	Let $\x:U \to \R^3$ be a proper wavefront, $\Omegam$ a tmb of $\x$, $\p \in \Sigma(\x)$ with $rank(D\x(\p))=0$ and $\Sigma(\x)_0=\{\q\in\Sigma(\x): rank(D\x(\q))=0\}$. If $\x$ is parallelly smoothable at $\p$ then:
	\begin{enumerate}
		
		\item There exist a open neighborhood $V$ of $\p$ in which one of the functions $k_1$, $k_2$ has a $C^\infty$ extension to $V-\Sigma(\x)_0$. More precisely, $k_1$ (resp. $k_2$) has a $C^\infty$ extension to $V-\Sigma(\x)_0$ if only if $\HO\leq0$ (resp. $\HO\geq0$) on $V$. 
		\item There exist a open neighborhood $V$ of $\p$ in which one of the functions $k_1$, $k_2$ converge to $\pm\infty$ (just one sign globally) near the singularities to $\p$. More precisely, $\lim\limits_{(u,v)\to \Sigma(\x)\cap V}k_1=\pm\infty$ (resp. $k_2$) if and only if $\HO\leq0$ (resp. $\HO\geq0$) on $V$.
	\end{enumerate}
\end{corollary}
\begin{proof}
	\
\begin{enumerate}

	\item Since that $\HO$ does not change sign on a neighborhood $V$ of $\p$ by item (2) of proposition \ref{p2} and applying proposition \ref{prin1} we get the result.
	\item By items (2) and (4) of proposition \ref{p2} $\HO$ and $k_1$ do not change sign on a neighborhood $V$ of $\p$ and applying proposition \ref{prin1} we get the result.  
\end{enumerate}
\end{proof}

\begin{proposition}
	Let $\x:U \to \R^3$ be a proper wavefront, $\Omegam$ a tmb of $\x$ and $\p \in \Sigma(\x)$ with $rank(D\x(\p))=0$. If $H$ is bounded on a neighborhood of $\p$  then we have:
	\begin{enumerate}
		\item There exist a neighborhood $V$ of $\p$ such that $K<0$ on $V-\Sigma(\x)$ and $\lim\limits_{(u,v)\to p}K=-\infty$.
		\item $\laO$ does not change sign on a neighborhood $V$ of $\p$.
		\item One of the function $k_1, k_2$ converge to $\infty$ and the another one to $-\infty$ at $\p$. More precisely, if $\laO\geq0$ (resp. $\laO\leq0$) on a neighborhood $V$ of $\p$ then $\lim\limits_{(u,v)\to p}k_1=-\infty$ (resp. $\infty$) and $\lim\limits_{(u,v)\to p}k_2=\infty$ (resp. $-\infty$). 
		\item There is no singularities of rank 1 on a neighborhood of $\p$.
		\item There exist a neighborhood $V$ of $\p$ such that $\lim\limits_{(u,v)\to \Sigma(\x)\cap V}\frac{k_1}{k_2}=-1$ 
	\end{enumerate}
\end{proposition}
\begin{proof}
	\
	\begin{enumerate}
		\item We know that $K\neq0$ near to $\p$. If there exist a sequence $\mathbf{a}_n\longrightarrow\p$ with $K(\mathbf{a}_n)>0$ for every $n\in \mathbb{N}$, as $H^2\geq K$ and $\lim\limits_{(u,v)\to p}K(\mathbf{a}_n)=\infty$ we have that $\lim\limits_{(u,v)\to p}|H(\mathbf{a}_n)|=\infty$ witch is contradictory, then we have (1). 
		\item By (1) $\laO\KO=\laO^2K<0$ near to $\p$ and since that $\KO(\p)\neq0$, then $\laO$ does not change sign on a neighborhood $V$ of $\p$.
		\item If $\laO\geq0$ near to $p$, $k_1=H-\sqrt{H^2-K}$ and $k_2=H+\sqrt{H^2-K}$ on a neighborhood of $\p$ and since that $K<0$ near to $\p$, then $k_2>0>k_1$ on a neighborhood of $\p$. By (3) of proposition \ref{propc2} we have the result.
		\item If $H$ is bounded on a neighborhood $V$ of $\p$ and suppose that exist a singularity $\q$ of rank $1$ in $V$, by (4) of proposition \ref{prin1} $\lim\limits_{(u,v)\to p}|H|=\infty$ witch is contradictory.
		\item Let $V$ be a bounded neighborhood of $\p$ with just singularities of rank 0 with $H$ bounded. There exist $C>0$ such that $|k_1+k_2|<C$, then $|1+\frac{k_1}{k_2}|<\frac{C}{|k_2|}$ and by (3) of proposition \ref{propc2} $\lim\limits_{(u,v)\to \q}\frac{k_1}{k_2}=-1$ for every $\q \in \Sigma(\x)\cap V$. Since that $\Sigma(\x)\cap V$ is compact we have the result.      
	\end{enumerate}
\end{proof}

\begin{proposition}\label{H1}
	Let $\x:U \to \R^3$ be a proper wavefront, $\Omegam$ a tmb of $\x$, $\p \in \Sigma(\x)$ with $rank(D\x(\p))=0$ and let us choose $W$ a compact neighborhood of $\p$ in which $\KO\neq 0$. The following statements are equivalents: 
	\begin{enumerate}
		\item The Mean curvature $H$ is bounded on $W-\Sigma(\x)$.
		\item There exist $C>0$ such that $|\HO|\leq C|\laO|$ on $W$. 
		\item There exist $C>0$ such that $|eG+gE-2fF|\leq C|\laO^2|$ on $W$.
		\item There exist $C>0$ such that $|\frac{1}{k_1}+\frac{1}{k_2}|\leq C|\laO|$ on $W$. 
		\item $\frac{1}{k_{1\Omega}}+\frac{1}{k_{2\Omega}}$ is bounded on $W$.  
	\end{enumerate}
	
\end{proposition}
\begin{proof}
	\
	\begin{itemize}
		\item $(1\Leftrightarrow2)$ Using that $\HO=\laO H$ on $W-\Sigma(\x)$ which is dense in $W$, we get the equivalence.
		\item $(1\Leftrightarrow3)$ Using that $H(EG-F^2)=eG+gE-2fF$ on $W-\Sigma(\x)$ and $EG-F^2\in \mathfrak{T}^2_\Omega(W)$ (see proposition \ref{D}), by compactness of $W$ we get the equivalence.
		\item $(2\Leftrightarrow4)$ by proposition \ref{propc2} $\frac{1}{k_1}+\frac{1}{k_2}$ has a $C^\infty$-extension to $W$ and this is equal to $\frac{\HO}{\KO}$. From this equality follows the equivalence. 
		\item $(4\Leftrightarrow5)$ Since that $k_{1\Omega}=\laO k_1$ and $k_{2\Omega}=\laO k_2$ on $W-\Sigma(\x)$ which is dense in $W$, we get the equivalence. 
	\end{itemize}
\end{proof}

\begin{proposition}
	Let $\x:U \to \R^3$ be a proper wavefront, $\Omegam$ a tmb of $\x$, $\p \in \Sigma(\x)$ with $rank(D\x(\p))=0$ and let us choose $V$ an open neighborhood of $\p$ in which $\KO\neq 0$. The following statements are equivalents: 
	\begin{enumerate}
	\item The Mean curvature $H$ has a $C^\infty$-extension to the neighborhood $V$ of $\p$.
	\item $\HO \in \mathfrak{T}_\Omega(V)$ 
	\item $eG+gE-2fF \in \mathfrak{T}^2_\Omega(V)$.
	\item $\frac{1}{k_1}+\frac{1}{k_2} \in \mathfrak{T}_\Omega(V)$ 
	\item $\frac{1}{k_{1\Omega}}+\frac{1}{k_{2\Omega}}$ has a $C^\infty$-extension to the neighborhood $V$ of $\p$.  
	\end{enumerate}
	 
\end{proposition}
\begin{proof}
	The proof of proposition \ref{H1} can be reproduced here to prove the corresponding equivalences.
\end{proof}

\begin{ex}\label{e1}
	The wavefront $\x:=(\frac{1}{2} \log \left(v^2+1\right)-\frac{1}{2} \log \left(u^2+1\right),\frac{u v}{v^2+1},\frac{u v^2}{v^2+1}-u+\tan ^{-1}(u))$ has as tmb:
	$$
	\Omegam=\begin{pmatrix}
	1&0\\
	0&1\\
	u&v
	\end{pmatrix}, \ \Lambdam=\begin{pmatrix}
	\frac{-u}{1+u^2}&\frac{v}{1+v^2}\\
	\frac{v}{1+v^2}&\frac{u(1-v^2)}{(1+v^2)^2}
	\end{pmatrix}, \ \mum=\begin{pmatrix}
	\frac{-(1+v^2)}{(1+u^2+v^2)^{\frac{3}{2}}}&\frac{uv}{(1+u^2+v^2)^{\frac{3}{2}}}\\
	\frac{uv}{(1+u^2+v^2)^{\frac{3}{2}}}&\frac{-(1+u^2)}{(1+u^2+v^2)^{\frac{3}{2}}}
	\end{pmatrix},$$
	then $\laO=\frac{-(u^2+v^2)}{(1+u^2)(1+v^2)^2}$, $\HO=-\frac{1}{2}(\lambda_{22}\mu_{11}-\lambda_{21}\mu_{12}+\lambda_{11}\mu_{22}-\lambda_{12}\mu_{21})=0$ and therefore the Mean curvature is extendable, with $H=0$ on $\R^2$.
\end{ex}

\section*{Acknowledgement}
I am grateful to  professors Kentaro Saji and Keisuke Teramoto for the
fruitful discussions about the extendibility of the Gaussian Curvature, during their visit to the ICMC-USP on September 2019 in which they gave the idea of expressing the extendibility property in terms of the ideal generated by $\laO$ and the second fundamental form.

\def\cprime{$'$}

\end{document}